\documentclass[10pt]{article}
\usepackage{amsmath,amssymb,latexsym,epsfig,subfigure,multirow,setspace}
\newtheorem{lem}{Lemma}[section]

\newtheorem{prop}[lem]{Proposition}

\newtheorem{thm}[lem]{Theorem}

\newtheorem{definition}[lem]{Definition}

\newtheorem{example}{Example}

\newcommand\pf{{\noindent\bf Proof.}~}
\newcommand\qed{\hfill\framebox[2.5mm]{}}

\newcommand{\Comment}[1]{}

\begin{document}

\title{Combinatorial Interpretation of General Eulerian Numbers}

\author{Tingyao Xiong\thanks{Radford University} \\ {\tt txiong@radford.edu}\\\\Jonathan I. Hall\thanks{Michigan State University}\\{\tt jhall@math.msu.edu}\\\\Hung-ping Tsao\thanks{Novato, CA} \\ {\tt tsaohp.tsao6@gmail.com }   }
\date{}

\maketitle

{\begin{abstract} \textbf{Since 1950s, mathematicians have successfully  interpreted  the traditional Eulerian numbers and $q-$Eulerian numbers combinatorially. In this paper, the authors give a combinatorial interpretation to the general Eulerian numbers defined on general arithmetic progressions $\{ a, a+d, a+2d,\dots\}$.}
\end{abstract}}


\section{Introduction}

\begin{definition}\label{excedance}
 Given a positive integer $n$, define $\Omega_n$ as the set of all permutations of $[n]=\{1,\; 2,\; 3,\;,\ldots, n\}$. For a permutation $\pi=p_1p_2p_3\ldots p_n\in \Omega_n$ , $i$ is called an ascent of $\pi$ if $p_i<p_{i+1}$; $i$ is called a weak excedance of $\pi$ if $p_i\ge i$.
\end{definition}

  It is well known that a traditional Eulerian number $A_{n, k}$ is the number of permtations $\pi\in \Omega_n$ that have $k$ weak excedances (\cite{Riordan}, page 215). And $A_{n,\;k}$ satisfies the recurrence:\;\;\;$A_{n,1}=1,\;(n\geq 1),\;\;A_{n,k}=0\;\;(k> n)$
\begin{equation}\label{recurrence1}
A_{n,k}=kA_{n-1,k}+(n+1-k)A_{n-1,k-1}\;\;(1\leq k\leq n)
\end{equation}

Besides the recursive formula (\ref{recurrence1}), $A_{n,\;k}$ can be calculated directly by the following analytic formula (\cite{Bona}, page 8):
\begin{equation}\label{analytic formula for Ank}
A_{n,k}=\sum_{i=0}^{k-1}(-1)^i(k-i)^n\left(\begin{array}{c}
                                                                 n+1 \\
                                                                 i
                                                               \end{array}
  \right)\hspace{1cm}(1\le k\le n)
\end{equation}
\begin{definition}\label{q eulerian}
Given a permutation $\pi=p_1p_2p_3\ldots p_n\in \Omega_n$, define functions
\[
maj\;\pi=\sum_{p_j>p_{j+1}}j\hspace{0.6cm}and
\]
\[
a(n,k,i)=\#\{\pi\;|\;maj\;\pi=i\;\; \&\;\;\pi\;\; has\;\;k\;\;ascents \}
\]
\end{definition}

Since the 1950's, Carlitz (\cite{Carlitz1},\cite{Carlitz2}) and his successors have generalized Euler's results to $q-$ sequences $\{1, q, q^2, q^3,\ldots, \}$. Under Carlitz's definition, the $q-$Eulerian numbers $A_{n,k}(q)$ are given by
\begin{equation}\label{q combinatorial1}
A_{n,k}(q)=q^{\frac{(m-k+1)(m-k)}{2}}\sum_{i=0}^{k(n-k-1)}a(n,n-k,i)q^i
\end{equation}
where functions $a(n,k,i)$ are as defined in Definition \ref{q eulerian}.

In \cite{XTH}, instead of studying $q-$sequences, the authors have generalized  Eulerian numbers to any general arithmetic progression
\begin{equation}\label{general arith}
\{a,\; a+d,\;a+2d,\;a+3d,\;\ldots \}
\end{equation}

Under the new definition, and given an arithmetic progression as defined in (\ref{general arith}), the general Eulerian numbers $A_{n,k}(a,d)$ can be calculated directly by the following equation (\cite{XTH}, Lemma 2.6)
\begin{equation}\label{formula of general eulerian numbers}
A_{n,k}(a,d)=\sum_{i=0}^{k}(-1)^i[(k+1-i)d-a]^n\left(\begin{array}{c}
n+1 \\
i
\end{array}
\right)
\end{equation}
Interested readers can find more results about the general Eulerian numbers and even general Eulerian polynomials in \cite{XTH}.

\section{Combinatorial Interpretation of General Eulerian Numbers}
The following concepts and properties will be heavily used in this section.

\begin{definition}\label{position of n}
Let $W_{n,k}$ be the set of $n-$permutations with $k$ weak excedances. Then $|W_{n,k}|=A_{n,k}$. Furthermore, given a permutation $\pi=p_1p_2p_3\ldots p_n$, let $Q_n(\pi)=i$ where $p_i=n$.
\end{definition}

Given a permutation $\pi\in \Omega_n$, it is known that $\pi$ can be written as a one line form like $\pi=p_1p_2p_3\ldots p_n$. Or $\pi$ can be written in a disjoint union of distinct cycles.  For $\pi$ written in a cycle form, we can use a \emph{standard representation} by writing (a) each cycle starting with its largest element, and (b) the cycles are written in increasing order of their largest element. Moreover, given a permutation $\pi$ written in a standard representation cycle form, define a function $f$ as $f(\pi)$ to be the permutation obtained from $\pi$ by erasing the parentheses. Then $f$ is known as the  \emph{fundamental bijection} from $\Omega_n$ to itself (\cite{Stanley}, page 30). Indeed, the inverse map $f^{-1}$ of the fundamental bijection function $f$ is also famous in illustrating the relation between the ascents and weak excedances as following: (\cite{Bona}, page 98)

\begin{prop}\label{ascents and weak excedances}
The function $f^{-1}$ gives a bijection between the set of permutations on $[n]$ with $k$ ascents and the set $W_{n,k+1}$.
\end{prop}

\begin{example}
The standard representation of permutation $\pi =5243716$ is $(2)(43)(7615)\in \Omega_7$, and $f(\pi)=2437615$; $Q_7(\pi)=5$; $\pi =5243716$ has $3$ ascents, while $f^{-1}(\pi)=(5243)(716)=6453271\in W_{7,4}$ has $3+1=4$ weak excedances because $p_1=6>1$, $p_2=4>2$, $p_3=5>3$, and $p_6=7>6$.
\end{example}

Now suppose we want to construct a sequence consisting of $k$ vertical bars and the first $n$ positive  integers. Then the $k$ vertical bars divide these $n$ numbers into $k+1$ compartments. In each compartment, there is either no number or all the numbers are listed in a decreasing order. The following definition is analogous to the definition of \cite{Bona}, page 8.

\begin{definition}
A bar in the above construction is called \emph{extraneous} if either
\begin{itemize}
  \item [(a)]it is immediately followed by another bar; or
  \item [(b)]after removing it each of the rest compartment either is empty or consists of integers in a decreasing order.
\end{itemize}
\end{definition}

\begin{example}
Suppose $n=7$, $k=4$, then in the following arrangement
\[
32\; |\; 1\; ||\; 7654\; |
\]
the $1$st, $2$nd, and the $4$th bars are extraneous.
\end{example}

Now we are ready to give combinatorial interpretations to the general Eulerian numbers $A_{n,k}(a,d)$. First note that equation (\ref{formula of general eulerian numbers}) implies that $A_{n,k}(a,d)$ is a homogeneous polynomial of degree $n$ with respect to $a$ and $d$. Indeed,
\begin{align}\label{Ank in the form of a and d-a}\nonumber
A_{n,k}(a,d)&=\sum_{i=0}^{k}(-1)^i[(k+1-i)d-a]^n\left(\begin{array}{c}
n+1 \\
i
\end{array}
\right)\\\nonumber
&=\sum_{i=0}^{k}(-1)^i[\;(k+1-i)(d-a)+(k-i)a\;]^n\left(\begin{array}{c}
n+1 \\
i
\end{array}
\right)\\\nonumber
&=\sum_{j=0}^n\left[ \sum_{i=0}^{k}(-1)^i(k+1-i)^{n-j}(k-i)^j\left(\begin{array}{c}
n+1 \\
i
\end{array}
\right)\right]\left(\begin{array}{c}
n \\
j
\end{array}
\right)(d-a)^{n-j}a^j\\
&=\sum_{j=0}^nc_{n,k}(j)\left(\begin{array}{c}
n \\
j
\end{array}
\right)(d-a)^{n-j}a^j
\end{align}
where
\begin{equation}\label{cnk formula}
c_{n,k}(j)=\sum_{i=0}^{k}(-1)^i(k+1-i)^{n-j}(k-i)^j\left(\begin{array}{c}
n+1 \\
i
\end{array}
\right), 0\le j\le n.
\end{equation}

The following Theorem gives combinatorial interpretations to the coefficients $c_{n,k}(j)$, $0\le j\le n$.
\begin{thm}\label{main theorem}
Let the general Eulerian numbers $A_{n,k}(a,d)$ be written as in equation $(\ref{Ank in the form of a and d-a})$. Then
\begin{equation}\label{cnk inter}
c_{n,k}(j)=\#\{\pi\in W_{n,k+1}\text{ and } j<Q_n(\pi)\le n\}+\#\{\pi\in W_{n,k}\text{ and } 1\le Q_n(\pi)\le j\}
\end{equation}
\end{thm}
\pf We can check the result in (\ref{cnk inter}) for two special values $j=0$ and $j=n$ quickly. By equation (\ref{analytic formula for Ank}),\\
\begin{itemize}
  \item []when $j=0$, $c_{n,k}(0)=\sum_{i=0}^{k}(-1)^i(k+1-i)^{n}\left(\begin{array}{c}
n+1 \\
i
\end{array}
\right)=A_{n,k+1}$;
\item []when $j=n$, $c_{n,k}(n)=\sum_{i=0}^{k}(-1)^i(k-i)^{n}\left(\begin{array}{c}
n+1 \\
i
\end{array}
\right)=A_{n,k}$. Therefore, (\ref{cnk inter}) is true for $j=0$ and $j=n$.
\end{itemize}

Generally, for $1\le j\le n-1$, we write down $k$ bars with $k+1$ compartments in between. Place each element of $[n]$ in a compartment. If none of the $k$ bars is extraneous, then the arrangement corresponds to a permutation with $k$ ascents. Let $B$ be the set of arrangements with at most one extraneous bar at the end and none of integers $\{1,2,\ldots,j\}$ locating in the last compartment. We will show that $c_{n,k}(j)=|B|$.\\

To achieve that goal, we use the Principle of Inclusion and Exclusion. There are $(k+1)^{n-j}k^j$ ways to put $n$ numbers into $k+1$ compartments with elements $\{1,2,\ldots,j\}$ avoiding the last compartments. Let $B_i$ be the number of arrangements with (1) none of $\{1,2,\ldots,j\}$ sitting in the last compartment; (2) there are at least $i$ ``separating" extraneous bars. Two extraneous bars are separating means that the two bars are not next to each other. Then the Principle of Inclusion and Exclusion shows that
\begin{equation}\label{B form}
|B|=(k+1)^{n-j}k^j-B_1+B_2+\ldots+(-1)^kB_k
\end{equation}

Now we consider the value of $B_i$, where $1\le i\le k$. Suppose that we have $k+1-i$ compartments with $k-i$ bars in between. There are $(k+1-i)^{n-j}(k-i)^j$ ways to insert $n$ numbers into these $k+1-i$ compartments with first $j$ integers avoiding the last compartment, and list integers in each component in a decreasing order. Then insert $i$ separating extraneous bars into $n+1$ positions. So we get
\begin{equation}\label{Bi form}
B_i=(k+1-i)^{n-j}(k-i)^j\left(\begin{array}{c}
n+1 \\
i
\end{array}
\right)
\end{equation}
Plug the formula (\ref{Bi form}) into equation (\ref{B form}), we have $c_{n,k}(j)=|B|$.

Given an arrangement $\pi\in B$, if we remove the bars, then we obtain a permutation $\pi\in\Omega_n$. So without confusion, we just use the same notation $\pi$ to represent both an arrangement in set $B$ and a permutation on $[n]$. Now for each $\pi\in B$, $\pi$  either
\begin{itemize}
  \item [case 1] has no extraneous bar and none of $\{1,2,\ldots,j\}$ locates in the last compartment; Or
  \item [case 2] has only one extraneous bar at the end.
\end{itemize}

If $\pi$ is in case 1, then $\pi$ has $k$ ascents since each bar is non-extraneous. And the last compartment of $\pi$ is nonempty. Therefore the last cycle of $f^{-1}(\pi)$ has to be $(n\ldots p_g)$. In other words, $Q_n(f^{-1}(\pi))=p_g>j$ since none of $\{1,2,\ldots, j\}$ locates in the last compartment. And by Proposition \ref{ascents and weak excedances}, $f^{-1}(\pi)\in W_{n,k+1}$.\\

If $\pi$ is in case 2, then $\pi$ has $k-1$ ascents since only the last bar is extraneous. Note that in this case, the arrangement with no elements of $\{1,2,\ldots,j\}$ in the compartment second to the last, or the last non-empty compartment have been removed by the Principle of Inclusion and Exclusion. Equivalently, at least one number of $\{1,2,\ldots,j\}$ has to be in the compartment second to the last. So the last cycle of $f^{-1}(\pi)$ has to be $(n\ldots p_l)$, and $Q_n(f^{-1}(\pi))=p_l\le j$. Also by Proposition \ref{ascents and weak excedances}, $f^{-1}(\pi)\in W_{n,k}$.

Combing all the results above, statement (\ref{cnk inter}) is correct. \qed
\vspace{0.1cm}
\newline

The next Theorem describes some interesting properties of the coefficients $c_{n,k}$.
\begin{thm}\label{second theorem}
Let the coefficients $c_{n,k}$ be as described in Theorem \ref{main theorem}. Then
\begin{itemize}
  \item [1.]$\sum_{k=0}^nc_{n,k}(j)=n!$, for any $0\le j\le n$;
  \item [2.]$c_{n,k}(j)=c_{n,n-k}(n-j)$, for all $0\le j,k\le n$
\end{itemize}
\end{thm}

Before we can prove Theorem \ref{second theorem}, we need the following Lemma which is also interesting by itself.
\begin{lem}\label{k and n+1-k}
Given a positive integer $n$, then
\[
\#\{\pi\in W_{n,k}\;\&\;Q_n(\pi)=j\}=\#\{\pi\in W_{n,n+1-k}\;\&\;Q_n(\pi)=n+1-j\}
\]
for any $1\le k, j\le n$.
\end{lem}
\pf First of all, given a positive integer $n$, we define a function $g: \Omega_n\rightarrow \Omega_n$ as following:
\[
\text{for }\pi=p_1p_2\ldots p_n\in\Omega_n,\;\; g(\pi)=(n+1-p_1)(n+1-p_2)\ldots (n+1-p_n)
\]
For instance, for $\pi=53214\in \Omega_5$, $g(\pi)=13452$. $g$ is obviously a bijection of $\Omega_n$ to itself.

Now for some fixed $1\le k, j\le n$, suppose $S_{k,j}=\{\pi\in W_{n,k}\;\&\;Q_n(\pi)=j\}$, and $T_{k,j}=\{\pi\in W_{n,n+1-k}\;\&\;Q_n(\pi)=n+1-j\}$. For any $\pi\in S_{k,j}$ written in the standard representation cycle form. So $\pi=(p_u\ldots)\ldots(n\ldots j)$ and $f(\pi)=p_u\ldots n\ldots j$ has $k-1$ ascents by Proposition \ref{ascents and weak excedances}. Now we compose $f(\pi)$ with the bijection function $g$ as just defined. Then $g(f(\pi))=n+1-p_u\ldots 1\ldots n+1-j$ has $n-k$ ascents, which implies that $f^{-1}(g(f(\pi)))$ has $n+1-k$ weak excedances. So $f^{-1}(g(f(\pi)))\in W_{n,n+1-k}$. Note that the last cycle of $f^{-1}(g(f(\pi)))$ has to be $(n\ldots n+1-j)$. Therefore $f^{-1}(g(f(\pi)))\in T_{k,j}$. Since both $f$ and $g$ are bijection functions, $f^{-1}gf$ gives a bijection between $S_{k,j}$ and $T_{k,j}$. \qed\newline

Now we are ready to prove Theorem \ref{second theorem}.\\

\pf (of Theorem \ref{second theorem})For part $1$, by Theorem \ref{main theorem},
\begin{align*}
\sum_{k=0}^nc_{n,k}(j)&=\sum_{k=0}^n\#\{\pi\in W_{n,k+1}\text{ and } j<Q_n(\pi)\le n\}+\sum_{k=0}^n\#\{\pi\in W_{n,k}\text{ and } 1\le Q_n(\pi)\le j\}\\
&=\sum_{k=0}^n \#\{\pi\in W_{n,k}\}=|\Omega_n|=n!
\end{align*}

For part $2$, also by Theorem \ref{main theorem},
\begin{align*}
c_{n,k}(j)&=\sum_{i=j+1}^n\#\{\pi\in W_{n,k+1}\text{ and } Q_n(\pi)=i\}+\sum_{m=1}^j\#\{\pi\in W_{n,k}\text{ and } Q_n(\pi)=m\}\\
&=\sum_{i=j+1}^n\#\{\pi\in W_{n,n-k}\text{ and } Q_n(\pi)=n+1-i\}\\
&+\sum_{m=1}^j\#\{\pi\in W_{n,n+1-k}\text{ and } Q_n(\pi)=n+1-m\}\hspace{3cm} \text {by Lemma }\ref{k and n+1-k}\\
&=\#\{\pi\in W_{n,k}\text{ and } 1\le Q_n(\pi)\le n-j\}+\#\{\pi\in W_{n,n+1-k}\text{ and } n-j< Q_n(\pi)\le n\}\\
&=c_{n,n-k}(n-j)\hspace{12cm}\qed
\end{align*}

\textbf{\emph{Remark}}\hspace{.5cm}Using the analytic formula of $c_{n,k}(j)$ as in (\ref{cnk formula}), part 2 of Theorem \ref{second theorem} implies the following identity:
\[
\sum_{i=0}^{k}(-1)^i(k+1-i)^{n-j}(k-i)^j\left(\begin{array}{c}
n+1 \\
i
\end{array}
\right)=\sum_{l=0}^{n-k}(-1)^l(n+1-k-l)^{j}(n-k-l)^{n-j}\left(\begin{array}{c}
n+1 \\
l
\end{array}
\right),
\]
where $n$ is a positive integer, and $0\le j,k\le n$.
\section{Another Combinatorial Interpretation of $c_{n,k}(1)$ and $c_{n,k}(n-1)$}
In pursuing the combinatorial meanings of the coefficients $c_{n,k}$, the authors have found some other interesting properties about permutations. The results in this section will reveal close connections between the traditional Eulerian numbers $A_{n,k}$ and $c_{n,k}(j)$, where $j=1$ or $j=n-1$.

One fundamental concept of permutation combinatorics is \emph{inversion}. A pair $(p_i,p_j)$ is called an \emph{inversion} of the permutation $\pi=p_1p_2\ldots p_n$ if $i<j$ and $p_i>p_j$ (\cite{Stanley}, page 36). The following definition provides the main concepts of this section.
\begin{definition}\label{AW and BW}
For a fixed positive integer $n$, let $AW_{n,k}=\{\pi=p_1p_2p_3\ldots p_n\;|\;\pi\in W_{n,k} \text{ and }p_1<p_n\}$ (or $(p_1, p_n)$ is not an inversion), and $BW_{n,k}=W_{n,k}\backslash AW_{n,k}$ (or $(p_1, p_n)$ is an inversion).
\end{definition}
It is obvious that $|AW_{n,k}|+|BW_{n,k}|=A_{n,k}$. The following Theorem interprets  coefficients  $c_{n,k}(1)$ and $c_{n,k}(n-1)$ in terms of $AW_{n,k}$ and $BW_{n,k}$.
\begin{thm}\label{cnk and aw, bw}
Let the coefficients $c_{n,k}$ of the general Eulerian numbers  be written as in equation $(\ref{cnk formula})$. $AW_{n,k}$ and $BW_{n,k}$ are as defined in Definition \ref{AW and BW}. Then
\begin{itemize}
  \item [(1)]$c_{n,k}(1)=2|AW_{n,k+1}|$;
  \item [(2)]$c_{n,k}(n-1)=2|BW_{n,k}|$.
\end{itemize}
\end{thm}
\pf For part (1). By Theorem \ref{main theorem}, $c_{n,k}(1)=|S_1|+|S_2|$, where $S_1=\{\pi=p_1p_2\ldots p_n|\pi\in W_{n,k+1}\;\&\;p_1\not=n\}$, $S_2=\{\pi=p_1p_2\ldots p_n|\pi\in W_{n,k}\;\&\;p_1=n\}$. Given a permutation $\pi=p_1p_2\ldots p_n\in S_1$ and $p_n\not=n$, then both $p_1p_2\ldots p_n$ and $p_np_2\ldots p_1$ belong to $S_1$, so one of them has to be in $AW_{n,k+1}$; If $\pi=p_1p_2\ldots p_n\in S_1$ and $p_n=n$, then $\pi\in AW_{n,k+1}$, but $p_np_2\ldots p_1\in S_2$. Therefore, $\frac{1}{2}c_{n,k}(1)=|AW_{n,k+1}|$.

Part (2) can be proved using exactly the same method. So we leave it to the readers as an exercise.\qed\\

$|AW_{n,k}|$ and $|BW_{n,k}|$  are interesting combinatorial concepts by themselves. Note that generally speaking, $|AW_{n,k}|\not=|BW_{n,k}|$. Indeed, $|AW_{n,k}|=|BW_{n,n+1-k}|$.
\begin{thm}\label{AWnk and BWnk}
For any positive integer $n\ge 2$, the sets $AW_{n,k}$ and $BW_{n,k}$ are defined in Definition \ref{AW and BW}. Then $|AW_{n,k}|=|BW_{n,n+1-k}|$ for $1\le k\le n$.
\end{thm}
\pf It is an obvious result of part 2 of Theorem \ref{second theorem} and Theorem \ref{cnk and aw, bw}.\qed\\

Our last result of this paper is the following Theorem which reveals that  both $|AW_{n,k}|$ and $|BW_{n,k}|$ take exactly the same recursive formula as the traditional Eulerian numbers $A_{n,k}$ as shown in equation (\ref{recurrence1}).

\begin{thm}\label{recursive formulae for AW and BW}
For a fixed positive integer $n$, let $AW_{n,k}$ and $BW_{n,k}$ be as defined in Definition  \ref{AW and BW}, then \begin{align}
k|AW_{n-1,k}|+(n+1-k)|AW_{n-1,k-1}|&=|AW_{n,k}|\hspace{1cm}\text{ and }\label{recursive formula for AW}\\
k|BW_{n-1,k}|+(n+1-k)|BW_{n-1,k-1}|&=|BW_{n,k}|\label{recursive formula for BW}
\end{align}
\end{thm}
\pf A computational proof can be obtained straightforward by using equation (\ref{cnk formula}) and Theorem \ref{cnk and aw, bw}. But here we provide a proof in a flavor of combinatorics.

\emph{Idea of the proof}: For equation (\ref{recursive formula for AW}), given a permutation $A_1=p_1p_2p_3\ldots p_{n-1}\in AW_{n-1,k}$, for each position $i$ with $p_{i}\ge i$, we insert $n$ into a certain place of  $A_1$, such that the new permutation $A_1'$ is in $AW_{n,k}$. There are $k$  such positions, so we can get $k$  new permutations in $AW_{n,k}$.  Similarly, if $A_2=p_1p_2p_3\ldots p_{n-1}\in AW_{n-1,k-1}$, for each position  $i$ with $p_{i}< i$, and the position at the end of $A_2$, we insert $n$ into a specific position of  $A_2$ and the resulting new permutation $A_2'$ is in $AW_{n,k}$. There are $n+1-k$ such positions, so we can get $n+1-k$  new permutations in $AW_{n,k}$. We will show that all the permutations obtained from the above constructions are distinct, and they have exhausted all the permutations in $AW_{n,k}$\\[1ex]

For any fixed $A'=\pi_1\pi_2\pi_3\ldots \pi_{n}\in AW_{n,k}$, then $\pi_1<\pi_n$. We classify $A'$ into the following disjoint cases:
\begin{itemize}
  \item [case a.]$\pi_i=n$ with $i<n$. So $A'=\pi_1\pi_2\ldots \pi_{i-1}n\pi_{i+1}\ldots \pi_{n-1}\pi_{n}$.
  \begin{itemize}
    \item [a.1]$\pi_1<\pi_{n-1}$, and $\pi_n\ge i$;
    \item [a.2]$\pi_1<\pi_{n-1}$, and $\pi_n < i$;
    \item [a.3]$\pi_1>\pi_{n-1}$, $\pi_n<n-1$, and $\pi_n\ge i$;
    \item [a.4]$\pi_1>\pi_{n-1}$, $\pi_n<n-1$, and $\pi_n< i$;
    \item [a.5]$\pi_1>\pi_{n-1}$, and $\pi_n=n-1$;
  \end{itemize}
  \item [case b.] $\pi_n=n$. So $\pi_i=n-1$ for some $i<n$ and $A'=\pi_1\pi_2\ldots \pi_{i-1}n-1\ldots \pi_{n-1}n$.
    \begin{itemize}
    \item [b.1]$\pi_1<\pi_{n-1}$;
    \item [b.2]$\pi_{n-1}<\pi_{1}<n-1$, and $\pi_{n-1} \ge i$;
    \item [b.3]$\pi_{n-1}<\pi_{1}<n-1$, and $\pi_{n-1} < i$;
    \item [b.4]$\pi_{1}=n-1$.
  \end{itemize}
\end{itemize}
Based on the classifications listed above, we can construct a map  $f: \{AW_{n-1,k},AW_{n-1,k-1}\}\rightarrow AW_{n,k}$ by applying the idea of the proof we have illustrated at the beginning of the proof. To save space, the map $f$ is demonstrated in Table \ref{Atable}. From Table \ref{Atable} we can see that in each case, the positions of inserting $n$ are all different. So all the images obtained in a certain case are different. Since all the cases are disjoint, all the images $A'\in AW_{n,k}$ are distinct.
\begin{table}[h]
  \centering
\begin{tabular}[h]{|c|c|c|c|}
\hline
&&&\\
$A=p_1p_2\dots p_{n-1}$& Position $i$& Condition& $A'\in AW_{n,k}$\\\hline
&&&\\
& &  &$A'=p_1p_2\dots p_{i-1} n p_{i+1}\dots p_{n-1}p_i$,\\
& & $p_i>p_1$ &with $p_1<p_{n-1}$ and $p_1< p_i$.\\
&&&(Case $a.1$)\\
\cline{3-4}
&&&\\
&$1<i\le n-1$&&$A'=p_1p_2\dots p_{i-1} n p_{i+1}\dots p_{n-2}p_ip_{n-1}$,\\
&and&$p_i<p_1$ and&with $p_i<p_1<p_{n-1}<n-1$.\\
&$p_i\ge i$&&(Case $a.3$)\\
\cline{3-4}
&&&\\
$A\in AW_{n-1,k}$&&&$A'=p_{1}p_2\dots p_{i-1} n-1 p_{i+1}\dots p_{n-2}p_{i}n$,\\
&&$p_i<p_1$ and&with $p_i<p_1$ and $p_i\ge i$.\\
&&&(Case $b.2$)\\
\cline{2-4}
&&&\\
&&$p_i=n-1$ and&$A'=p_{n-1}p_2\dots p_{i-1}np_{i+1}\dots p_{n-2}p_{1}n-1$,\\
&&$i<n-1$ &with $p_1<p_{n-1}$.\\
&$i=1$&&(Case $a.5$)\\
\cline{3-4}
&&&\\
&&$p_{n-1}=n-1$ &$A'=n-1p_2\dots p_{n-2}p_1n$,\\
&&&(Case $b.4$)\\\hline
&&&\\
& &  &$A'=p_1p_2\dots p_{i-1} n p_{i+1}\dots p_{n-1}p_i$,\\
& & $p_i>p_1$ &with $p_1<p_{n-1}$ and $p_1<p_i$.\\
&&&(Case $a.2$)\\
\cline{3-4}
&&&\\
&$1<i\le n-1$&$p_i<p_1$ and&$A'=p_1p_2\dots p_{i-1} n p_{i+1}\dots p_{n-2}p_ip_{n-1}$,\\
&and&$p_{n-1}<n-1$&with $p_i<p_1<p_{n-1}<n-1$.\\
&$p_i<i$&&(Case $a.4$)\\
\cline{3-4}
&&&\\
$A\in AW_{n-1,k-1}$&&$p_i<p_1$ and&$A'=p_1p_2\dots p_{i-1}n-1p_{i+1}\dots p_{n-2} p_i n$.\\
&&$p_{n-1}=n-1$&with $p_i<p_1$ and $p_{i}<i$.\\
&&&(Case $b.3$)\\
\cline{2-4}
&&&\\
&$i=n$&&$A'=p_1p_2\dots p_{n-1}n$.\\
&& &(Case $b.1$)\\\hline
 \end{tabular}
 \caption{The Map $f: \{AW_{n-1,k},AW_{n-1,k-1}\}\rightarrow AW_{n,k}$}
 \label{Atable}
\end{table}

Similarly, for each $B'=\pi_1\pi_2\pi_3\ldots \pi_{n}\in BW_{n,k}$, then $\pi_1>\pi_n$. We classify $B'$ into the following disjoint cases:
\begin{itemize}
  \item [case c.]$\pi_i=n$ with $1<i\le n-1$. So $B'=\pi_1\pi_2\ldots \pi_{i-1}n\pi_{i+1}\ldots \pi_{n-1}\pi_n$.
  \begin{itemize}
    \item [c.1]$\pi_1>\pi_{n-1}$, and $\pi_n\ge i$;
    \item [c.2]$\pi_1>\pi_{n-1}$, and $\pi_n < i$;
    \item [c.3]$\pi_1<\pi_{n-1}<n-1$, $\pi_{n-1}\ge i$;
    \item [c.4]$\pi_1<\pi_{n-1}<n-1$, $\pi_{n-1}< i$;
    \item [c.5]$\pi_{n-1}=n-1$;
    \item [c.6]$\pi_{n-1}=n$;
  \end{itemize}
  \item [case d.] $\pi_1=n$. So  $B'=n\pi_2\ldots \pi_{n-2}\pi_{n-1}$.
    \begin{itemize}
    \item [d.1]$\pi_{n-2}<\pi_{n-1}$;
    \item [d.2]$\pi_{n-2}>\pi_{n-1}$.
  \end{itemize}
\end{itemize}
To prove equation (\ref{recursive formula for BW}), we use a similar idea of proof as shown above. If $B_1=p_1p_2p_3\ldots p_{n-1}\in BW_{n-1,k}$, for each position $i$ with $p_{i}\ge i$, we insert $n$ into a certain place of  $B_1$ to get $B_1'\in AW_{n,k}$; If $B_2=p_1p_2p_3\ldots p_{n-1}\in BW_{n-1,k-1}$, for each position  $i$ with $p_{i}< i$, and the position $i$ where $p_i=n-1$, we insert $n$ into a specific position of  $B_2$ to obtain $B_2'\in AW_{n,k}$. Such a map $g: \{BW_{n-1,k},BW_{n-1,k-1}\}\rightarrow BW_{n,k}$ is illustrated in Table \ref{Btable}. And the distinct images under $g$ exhaust all the permutations in $BW_{n,k}$.
\begin{table}[h]
  \centering
\begin{tabular}{|c|c|c|c|}
\hline
&&&\\
$B=p_1p_2\dots p_{n-1}$& Position $i$& Condition& $B'\in BW_{n,k}$\\\hline
&&&\\
& &  &$B'=p_1\dots p_{i-1} n p_{i+1}\dots p_{n-1}p_i$,\\
& & $p_1>p_i$ &with $p_1>p_{n-1}$ and $p_i\ge i$.\\
&&&(Case $c.1$)\\
\cline{3-4}
&&&\\
&$1<i<n-1$&&$B'=p_1\dots p_{i-1} n p_{i+1}\dots p_ip_{n-1}$,\\
&and&$p_1<p_i<n-1$&with $p_1<p_i<n-1$ and $p_i\ge i$.\\
&$p_i\ge i$&&(Case $c.3$)\\
\cline{3-4}
&&&\\
$B\in BW_{n-1,k}$&&&$B'=np_2\dots p_{i-1} n-1 p_{i+1}\dots p_{1}p_{n-1}$,\\
&&$p_1<p_i=n-1$&with $p_i=n-1$ and $p_1> p_{n-1}$.\\
&&&(Case $d.2$)\\
\cline{2-4}
&&&\\
&&&$B'=np_2\dots p_{n-1}p_1$,\\
&$i=1$&$p_1\ge 1$ &with $p_{n-1}<p_1$.\\
&&&(Case $d.1$)\\\hline
&&&\\
& &  &$B'=p_1\dots p_{i-1} n p_{i+1}\dots p_{n-1}p_i$,\\
& & $p_1>p_i$ &with $p_1>p_{n-1}$ and $p_i< i$.\\
&$1<i<n-1$&&(Case $c.2$)\\
\cline{3-4}
&and&&\\
&$p_i< i$&&$B'=p_1\dots p_{i-1} n p_{i+1}\dots p_ip_{n-1}$,\\
&&$p_1<p_i$&with $p_1<p_i<n-1$ and $p_i< i$.\\
&&&(Case $c.4$)\\
\cline{2-4}
&&&\\
$B\in BW_{n-1,k-1}$&$i=n-1$&$p_1>p_i=p_{n-1}$&$B'=p_1p_2\dots p_{n-2} n p_{n-1}$.\\
&$p_i<i$&&(Case $c.6$)\\
&&&\\
\cline{2-4}
&&&\\
&$1\le i<n-1$&&$B'=p_1p_2\dots p_{i-1}np_{i+1}\dots p_{n-2}n-1p_{n-1}$.\\
&and&$p_i=n-1$ &(Case $c.5$)\\
&$p_i=n-1$&&\\\hline
 \end{tabular}
 \caption{The Map  $g: \{BW_{n-1,k},BW_{n-1,k-1}\}\rightarrow BW_{n,k}$}
 \label{Btable}
\end{table}

Here is a concrete example for the constructions illustrated in Table \ref{Btable}:
\begin{example}
Suppose $n=4$, $k=2$. We want to obtain $BW_{4,2}=\{3142, 3412, 3421, 4132, 4213, 4312, 4321\}$ from $BW_{3,2}=\{321, 231\}$ and $BW_{3,1}=\{312\}$. For $321\in BW_{3,2}$, $p_1=3\ge 1$, then it corresponds to $B'=4213$ which is Case d.1 in Table \ref{Btable}; $p_2=2\ge 2$, then it corresponds to $B'=3412$ which is Case c.1 in Table \ref{Btable}. Similarly, we can construct $\{4312, 4321\}$ from $231\in BW_{3,2}$, and $\{3421, 3142,4132\}$ from $312\in BW_{3,1}$ using Table \ref{Btable}.
\end{example}


\end{document}